\documentclass[twoside,12pt]{article}
\usepackage{amsmath}
\usepackage{amssymb}
\usepackage[usenames,dvipsnames]{color}
\definecolor{rosso}{rgb}{0.8,0,0}

\title{{C}ahn--{H}illiard approach
to some degenerate parabolic equations 
with dynamic boundary conditions}

\author{Takeshi Fukao\\
Department of Mathematics, Faculty of Education\\
Kyoto University of Education\\
1~Fujinomori, Fukakusa, Fushimi-ku, Kyoto~612-8522 Japan\\
E-mail: \texttt{fukao@kyokyo-u.ac.jp}}
 
\date{}

\newcommand\testopari{\sc Takeshi Fukao}
\newcommand\testodispari{\sc {S}tefan problem with dynamic boundary conditions}
\begin{document}

\maketitle

\begin{abstract}
In this paper 
the well-posedness of some degenerate parabolic equations 
with a dynamic boundary condition is considered.
To characterize the target degenerate parabolic equation from 
the {C}ahn--{H}illiard system, 
the nonlinear term coming from the convex part of the double-well
potential is chosen using a suitable maximal monotone graph. 
The main topic of this paper is the existence problem under an 
assumption for this maximal monotone graph for treating a wider class. 
The existence of a weak solution is proved.

\vspace{2mm}

\vspace{2mm}
\noindent \textbf{AMS (MOS) subject clas\-si\-fi\-ca\-tion: 35K65, 35K30, 47J35} 

\end{abstract}

\section{Introduction}
\setcounter{equation}{0}

The relationship between the {A}llen--{C}ahn equation \cite{AC79} and the motion 
by mean curvature 
is interesting 
as the singular limit of the following form: 
\begin{equation} 
	\frac{\partial u}{\partial t}-\Delta u 
	+ \frac{1}{\varepsilon ^2}(u ^3-u)=0 
	\quad {\rm in}~Q:=(0,T) \times \Omega, 
	\nonumber 
\end{equation} 
as $\varepsilon \searrow 0$, where $0<T<+\infty$ 
and $\Omega \subset \mathbb{R}^d$ for $d=2,3$, which is a bounded domain with smooth 
boundary $\Gamma$. 
For example, {B}ronsard and {K}ohn presented a pioneering result in \cite{BK91}, 
and subsequently many related results have been obtained. 
A similar concept in this framework, the {C}ahn--{H}illiard system \cite{CH58},
is connected to 
motion by the {M}ullins--{S}ekerka law \cite{MS63} in the limit of
\begin{gather}
	\frac{\partial u}{\partial t}-\Delta \mu  =0 \quad {\rm in}~Q, 
	\nonumber \\
	\mu =- \varepsilon \Delta u  + \frac{1}{\varepsilon }(u ^3-u) \quad {\rm in}~Q
	\label{MS}
\end{gather}
as $\varepsilon \searrow 0$. 
For both of these, the target problems are sharp interface models in a classical sense and 
a powerful analysis tool seems to be the method of matched asymptotic expansions
(see\ \cite{ABC94, CC98, Peg89} and the references in these papers).

In this paper, we discuss this relation from a different view point. 
To do so, we begin with the following degenerate parabolic equation:
\begin{equation} 
	\frac{\partial u}{\partial t} -\Delta \beta (u) = g 
	\quad {\rm in}~Q, 
	\label{dpe}
\end{equation}
where $g:\Omega \to \mathbb{R}$ is a given source. 
This equation is characterized by the choice of $\beta: \mathbb{R} \to \mathbb{R}$. 
For example, if 
we choose $\beta$ to be a piecewise linear function 
of the form 
\begin{equation} 
	\beta (r):=
	\begin{cases}
	k_s r & r<0, \\
	0     & 0\le r \le L, \\
	k_\ell (r-L) & r >L;
	\end{cases} 
	\label{stefan}
\end{equation}
where $k_s$ and $k_\ell >0$ represent the heat conductivities of the 
solid and liquid regions, respectively, and $L>0$ is the latent heat constant, 
then \eqref{dpe} is the weak formulation of the {S}tefan problem, 
or the ``enthalpy formulation,'' 
where the unknown $u$ denotes the enthalpy and 
$\beta (u)$ denotes the temperature. 
The informant of the sharp interface, in other words the {S}tefan condition, is 
hidden in the weak formulation. 
Another example 
is the weak formulation of the {H}ele-{S}haw problem. 
If we choose $\beta $ to be the inverse of the {H}eaviside function 
\begin{equation*}
	{\mathcal H}(r):=
	\begin{cases}
	\displaystyle 0 & {\rm if~}r<0, \vspace{2mm}\\
	\displaystyle [0,1] & {\rm if~}r=0, \vspace{2mm}\\
	\displaystyle 1 & {\rm if~}r >0
	\end{cases} 
	\quad {\rm for~all~} r \in \mathbb{R},
\end{equation*}
so that $\beta$ is the multivalued function 
$\beta(r):={\mathcal H}^{-1}(r)=\partial I_{[0,1]}(r)$
for all $r \in [0,1]$,
then \eqref{dpe} can be stated as
\begin{equation*} 
	\xi  \in \beta (u), \quad 
	\frac{\partial u}{\partial t} -\Delta \xi = g
	\quad {\rm in}~Q,
\end{equation*}
where $\partial I_{[0,1]}$ is 
the subdifferential of the indicator function $I_{[0,1]}$ on the interval $[0,1]$, 
the unknown $u$ denotes the order parameter.
Details about weak formulations may be found in {V}isintin \cite{Vis96}. 
Weak formulations for this kind of sharp interface model are 
the focus of this paper. 
Therefore, we use the terms ``{S}tefan problem'' and ``{H}ele-{S}haw problem'' 
in the sense of weak formulations throughout this paper.

Recently, the author considered the approach to the following {C}ahn--{H}illiard system 
for the {S}tefan problem in \cite{Fuk15}: 
\begin{gather} 
	\frac{\partial u}{\partial t} -\Delta \mu = 0 
	\quad \mbox{in }Q, 
	\label{int1} \\
	\mu = -\varepsilon \Delta u + \beta (u) 
	+ \varepsilon \pi (u)-f 
	\quad \mbox{in }Q,
	\label{int2}
\end{gather}
with a dynamic boundary condition of the form
\begin{gather} 
	\frac{\partial u}{\partial t} + \partial_{\mbox{\scriptsize \boldmath $ \nu $}} \mu - \Delta _\Gamma \mu  =0
	\quad \mbox{on }\Sigma:=(0,T) \times \Gamma,
	\label{int3}
	\\
	\mu =\varepsilon\partial _{\mbox{\scriptsize \boldmath $ \nu $}} u - \varepsilon\Delta _\Gamma u 
	+ \beta (u)
	+ \varepsilon\pi (u)-f_\Gamma 
	\quad \mbox{on }\Sigma,
	\label{int4}
\end{gather}
where the symbol 
$\partial_{\mbox{\scriptsize \boldmath $ \nu $}} $ 
denotes the normal derivative on the boundary $\Gamma $ outward from $\Omega$,
the symbol $\Delta _{\Gamma }$ stands for the {L}aplace--{B}eltrami operator
on $\Gamma $ (see, e.g., \cite[Chapter~3]{Gri09}), 
$\beta $ is defined by \eqref{stefan}, and 
$\pi :\mathbb{R} \to \mathbb{R}$ is a piecewise linear function 
defined by $\pi (r):=L/2$ if $r<0$, 
$\pi (r):=L/2-r$ if $0 \le r \le L$ and 
$\pi (r):=-L/2$ if $r>L$.
Thanks to this choice, system \eqref{int1}--\eqref{int4} 
has the structure of a {C}ahn--{H}illiard system. 
This problem originally comes from \cite{GMS11}. 
Formally, if we let $\varepsilon \searrow 0$ in 
\eqref{int1}--\eqref{int4}, then we can see that 
the {C}ahn--{H}illiard system \eqref{int1}--\eqref{int4} converges in a suitable sense to the following 
{S}tefan problem with a dynamic boundary condition:
\begin{gather*} 
	\frac{\partial u}{\partial t} 
	-\Delta \beta (u) = -\Delta f 
	\quad \mbox{in }Q, 
	\\
	\frac{\partial u}{\partial t} 
	+ \partial_{\mbox{\scriptsize \boldmath $ \nu $}} \beta (u)
	- \Delta _\Gamma \beta (u)  =\partial_{\mbox{\scriptsize \boldmath $ \nu $}} f_\Gamma 
	- \Delta _\Gamma f_\Gamma 
	\quad \mbox{on }\Sigma.
\end{gather*}
Here, we should take care of the difference 
between the order and position of $\varepsilon $
in \eqref{MS} and 
\eqref{int2} even when $\beta(u)=u^3$ and $\pi (u)=-u$. 
In \cite{Fuk15}, $\beta $ is assumed to satisfy the following condition:

\begin{quote}
$\beta$ is a maximal monotone graph in $\mathbb{R}\times \mathbb{R}$, 
and is a subdifferential $\beta =\partial \widehat{\beta }$ 
of some proper, lower semicontinuous, and convex function 
$\widehat{\beta}:\mathbb{R} \to [0,+\infty ]$
satisfying $\widehat{\beta }(0)=0$ with some effective domain $D(\beta )$. This implies 
$\beta (0)=0$. Moreover, 
there exist two constants $c$, $\tilde{c}>0$ such that 
\begin{equation} 
	\widehat{\beta} (r) \ge c |r|^2- \tilde{c} \quad {\rm for~all}~r \in \mathbb{R}.
	\label{akey}
\end{equation}
\end{quote}
It is easy to see that \eqref{dpe} 
represents a large number of problems, 
including the 
porous media equation, the nonlinear diffusion equation of {P}enrose--{F}ife type, 
the fast diffusion equation, and so on. However, 
to apply this approach from the {C}ahn--{H}illiard system 
to these wider classes of the degenerate parabolic equation, the 
growth condition \eqref{akey} 
is too strong (see, e.g.\ \cite{DK99}). 
Therefore, in this paper based on the essential idea from 
\cite{CF16},  
we relax the assumption in \eqref{akey}. 
This is the different point from the previous work \cite{Fuk15}. See also 
\cite{AMTI06, Igb07, IK02} for related problems of interest.

\paragraph{Notation.}
Let $H:=L^2(\Omega )$, $V:=H^1(\Omega )$, 
$H_\Gamma :=L^2(\Gamma )$ and $V_\Gamma :=H^1(\Gamma )$
with the usual norms 
$| \cdot |_{H}$, $|\cdot |_{V}$, $|\cdot |_{H_\Gamma}$, $|\cdot |_{V_\Gamma}$ 
and inner products $(\cdot,\cdot )_{H}$, $(\cdot ,\cdot )_{V}$,
$(\cdot,\cdot )_{H_\Gamma}$, $(\cdot ,\cdot )_{V_\Gamma}$, respectively, and let 
$\boldsymbol{H}:=H \times H_\Gamma$, 
$\boldsymbol{V}:=\{ (z,z_\Gamma ) \in V \times V_\Gamma 
: z_\Gamma =z_{|_\Gamma} \mbox{ a.e.\ on }\Gamma  \}$ and 
$\boldsymbol{W}:=H^2(\Omega ) \times H^2(\Gamma )$. 
Then $\boldsymbol{H}$, $\boldsymbol{V}$ and $\boldsymbol{W}$ are 
{H}ilbert spaces with 
the inner product
\begin{equation*}
	(\boldsymbol{u},\mbox{\boldmath $ z $}
	)_{\boldsymbol{H}}
	:=(u,z)_{H} + (u_\Gamma ,z_\Gamma )_{H_\Gamma } \quad 
	\mbox{for all}~\boldsymbol{u}, 
	\boldsymbol{z}
	\in \boldsymbol{H},
\end{equation*}
and the related norm 
is analogously defined as one of 
$\boldsymbol{V}$ or $\boldsymbol{W}$. 
Define $m:\boldsymbol{H} \to \mathbb{R}$ by 
\begin{equation*}
	m(\boldsymbol{z}):=\frac{1}{|\Omega |+|\Gamma| }
	\left\{ \displaystyle \int_{\Omega }^{}z dx
	+ \int_{\Gamma }^{} z_{\Gamma } d\Gamma \right\} 
	\quad \mbox{for all }\boldsymbol{z} \in \boldsymbol{H},
\end{equation*}
where $|\Omega |:=\int_{\Omega }^{}1dx$ and $|\Gamma |:=\int_{\Gamma }^{}1d\Gamma $. 
The symbol $\boldsymbol{V}^*$ denotes 
the dual space of $\boldsymbol{V}$, 
and the pair 
$\langle \cdot ,\cdot 
\rangle _{\boldsymbol{V}^*,\boldsymbol{V}}$
denotes the duality pairing between $\boldsymbol{V}^*$ and 
$\boldsymbol{V}$. 
Moreover, define 
the bilinear form 
$a(\cdot ,\cdot ):\boldsymbol{V} \times \boldsymbol{V} \to \mathbb{R}$ by 
\begin{equation*}
	a(\boldsymbol{u},\boldsymbol{z}):=
	\int_{\Omega }^{} \nabla u \cdot \nabla z dx 
	+\int_{\Gamma }^{} \nabla _\Gamma u _\Gamma \cdot \nabla _\Gamma z_\Gamma d\Gamma 
	\quad \mbox{for all }\boldsymbol{u},
	\boldsymbol{z} \in \boldsymbol{V},
\end{equation*} 
where 
$\nabla _{\Gamma }$ denotes the surface gradient on 
$\Gamma$ (see, e.g., \cite[Chapter~3]{Gri09}). 
We introduce 
the subspace 
$\boldsymbol{H}_0:=\{ \boldsymbol{z} \in
\boldsymbol{H} : m(\boldsymbol{z})=0 \}$ of $\boldsymbol{H}$
and $\boldsymbol{V}_0 :=\boldsymbol{V} \cap \boldsymbol{H}_0$, with 
their norms $| \boldsymbol{z}
|_{\boldsymbol{H}_0}:=|\boldsymbol{z}
|_{\boldsymbol{H}}$ for all 
$\boldsymbol{z} \in \boldsymbol{H}_0$ and 
$|\boldsymbol{z}|_{\boldsymbol{V}_0}:=
a(\boldsymbol{z},\boldsymbol{z} )^{1/2}
$
for all 
$\boldsymbol{z} \in \boldsymbol{V}_0$. 
Then the duality mapping 
$\boldsymbol{F}: \boldsymbol{V}_0 \to \boldsymbol{V}_0^*$ is defined by 
$\langle \boldsymbol{F}
\boldsymbol{z}, \tilde{\boldsymbol{z}} 
\rangle _{\boldsymbol{V}_0^*, \boldsymbol{V}_0}
:= a(\boldsymbol{z},\tilde{\boldsymbol{z}})
$
for all 
$\boldsymbol{z}, \tilde{\boldsymbol{z}} \in \boldsymbol{V}_0$ 
and the inner product in $\boldsymbol{V}_0^*$ is defined by 
$
(\boldsymbol{z}_1^*,\boldsymbol{z}_2^*)_{\boldsymbol{V}_0^*}
:=\langle \boldsymbol{z}_1^*, 
\mbox{\boldmath $ F$} ^{-1} \boldsymbol{z}_2^* 
\rangle _{\boldsymbol{V}_0^*,\boldsymbol{V}_0}
$
for all 
$\boldsymbol{z}_1^*,\boldsymbol{z}_2^* \in \mbox{\boldmath $ V$}_0^*$. 
Moreover, define $\boldsymbol{P}:\boldsymbol{H} \to \boldsymbol{H}_0$ by 
$\boldsymbol{P}\boldsymbol{z}:=\boldsymbol{z}-m(\boldsymbol{z})\boldsymbol{1}$ 
for all $\boldsymbol{z} \in \boldsymbol{H}$, where $\boldsymbol{1}:=(1,1)$. 
Thus we obtain the dense and compact embeddings
$\boldsymbol{V}_0 
\mathop{\hookrightarrow} \mathop{\hookrightarrow}
\boldsymbol{H}_0 
\mathop{\hookrightarrow} \mathop{\hookrightarrow}
\boldsymbol{V}_0^*$. See \cite{CF14, CF15} for further details.

\section{Existence of the weak solution}
\setcounter{equation}{0}

In this section, 
we state an existence theorem 
for the weak solution of a degenerate parabolic equation 
with a dynamic boundary condition of the following form:
\begin{gather} 
	\xi \in \beta (u), \quad \frac{\partial u}{\partial t}-\Delta \xi  = g 
	\quad \mbox{\rm a.e.\ in }Q,
	\nonumber 
\\
	\xi _\Gamma \in \beta (u_\Gamma ), 
	\quad  
	\xi _\Gamma = \xi_{|_\Gamma },
	\quad 
	\frac{\partial u_\Gamma }{\partial t}
	+ \partial_{\mbox{\scriptsize \boldmath $ \nu $}} 
	\xi 
	- \Delta _\Gamma \xi _\Gamma =g_\Gamma 
	\quad \mbox{\rm a.e.\ on }\Sigma,
	\nonumber 
\\
	u(0)=u_0
	\quad 
	\mbox{\rm a.e.\ in }\Omega, \quad 
	u_\Gamma  (0)=u_{0\Gamma}
	\quad \mbox{\rm a.e.\ on }\Gamma,
	\nonumber 
\end{gather} 
where $\beta$, $g$, $g_\Gamma$, $u_0$ and $u_{0\Gamma }$ 
satisfy the following assumptions:
\begin{itemize}
 \item[(A1)] $\beta$ is a maximal monotone graph in $\mathbb{R}\times \mathbb{R}$, 
and is a subdifferential $\beta =\partial \widehat{\beta }$ 
of some proper, lower semicontinuous, and convex function 
$\widehat{\beta}:\mathbb{R} \to [0,+\infty ]$
satisfying $\widehat{\beta }(0)=0$ in some effective domain $D(\beta )$. This implies that
$\beta (0)=0$;
 \item[(A2)] $\boldsymbol{g}\in L^2(0,T;\boldsymbol{H}_0)$; 
 \item[(A3)] $\boldsymbol{u}_0:=(u_0,u_{0\Gamma}) \in \boldsymbol{H}$ with $m_0 \in {\rm int}D(\beta )$, and 
the compatibility conditions $\widehat{\beta }(u_0)\in L^1(\Omega ), 
\widehat{\beta }(u_{0\Gamma }) \in L^1(\Gamma)$ hold. 
\end{itemize}
We remark that the growth condition of $\widehat{\beta}$ in 
(A1) and the regularity of $\boldsymbol{u}_0$ in (A3) 
are relaxations from a previous related result \cite{Fuk15} (cf.\ \eqref{akey}).

\paragraph{Theorem 2.1.} {\it Under assumptions 
{\rm (A1)}--{\rm (A3)}, 
there exists at least one pair 
$(\boldsymbol{u}, \boldsymbol{\xi })$ of 
functions $\boldsymbol{u} \in H^1(0,T;\boldsymbol{V}^*) 
\cap L^2(0,T;\boldsymbol{H})$ and 
$\boldsymbol{\xi } \in L^2(0,T;\boldsymbol{V})$
such that 
$\xi \in \beta (u)$ a.e.\ in $Q$, $\xi _\Gamma \in \beta (u_\Gamma )$ and 
$\xi _\Gamma =\xi _{|_\Gamma }$ a.e.\ on $\Sigma$, and that satisfy}
\begin{align} 
	\bigl \langle u' (t),z
	\bigr \rangle _{V^*,V}
	& {}+ \bigl \langle u_\Gamma ' (t),z_\Gamma 
	\bigr \rangle _{V^*_\Gamma ,V_\Gamma }
	+ \int_{\Omega }^{} \nabla \xi (t) \cdot \nabla z dx 
	+ \int_{\Gamma }^{} \nabla_\Gamma  \xi _{\Gamma}(t) 
	\cdot \nabla _\Gamma z_\Gamma 
	d \Gamma 
	\nonumber \\
	& 
	= \int_{\Omega }^{} g(t) z dx 
	+ \int_{\Gamma }^{} g_\Gamma (t) z_\Gamma 
	d \Gamma 
	\quad {\it for~all~}\boldsymbol{z}:=(z,z_\Gamma ) \in \boldsymbol{V}
	\label{def}
\end{align}
{\it for a.a.\ $t \in (0,T)$ with 
$u(0)=u_0$
a.e.\ in $\Omega$ and 
$u_\Gamma  (0)=u_{0\Gamma}$
a.e.\ on $\Gamma$.} \\

The continuous dependence is completely the same as in a previous result 
\cite[Theorem~2.2]{Fuk15}. Therefore, we devolve the uniqueness problem  on 
\cite{Fuk15}.

\section{Proof of the main theorem}
\setcounter{equation}{0}

In this section, we prove the main theorem. 
The strategy of the proof is similar to that of \cite[Theorem~2.1]{Fuk15}. 
However, to relax the assumption we use 
a different uniform estimate. 
Let us start with an approximate problem. 
Recall the {Y}osida approximation 
$\beta _\lambda :\mathbb{R} \to \mathbb{R}$  
and the related {M}oreau--{Y}osida 
regularization $\widehat{\beta }_\lambda$
of $\widehat{\beta }:\mathbb{R} \to \mathbb{R}$
(see, e.g., \cite{Bar10}).
We see that 
$0 \le \widehat{\beta }_\lambda (r) 
\le \widehat{\beta }(r)$
for all $r \in \mathbb{R}$. 
Moreover, we define the following proper, 
lower semicontinuous, and convex functional 
$\varphi: \boldsymbol{H}_0 \to [0,+\infty ]$:
\begin{equation}
	\varphi (\boldsymbol{z}) 
	:= \begin{cases}
	\displaystyle 
	\frac{1}{2} \int_{\Omega }^{} | \nabla z |^2 dx 
	+\frac{1}{2} \int_{\Gamma }^{} 
	 |\nabla _{\Gamma } z_{\Gamma }  |^2 d\Gamma  
	\quad 
	\mbox{if } 
	\boldsymbol{z} \in \boldsymbol{V}_0, \vspace{2mm}\\
	+\infty \quad \mbox{otherwise}. 
	\end{cases} 
	\nonumber 
\end{equation}
The subdifferential 
$\partial \varphi $ on $\boldsymbol{H}_0$ is characterized by 
$\partial \varphi (\boldsymbol{z})
=(-\Delta z,\partial_{\boldsymbol{\nu }} z-\Delta _\Gamma z_\Gamma ) 
$ with 
$\boldsymbol{z} \in 
D(\partial \varphi )=\boldsymbol{W} \cap \mbox{\boldmath $ V$}_0$ 
(see, e.g., \cite[Lemma~C]{CF15}). 
By virtue of the well-known theory of evolution equations
(see, e.g., \cite{CF14, CF15, CV90, KNP95}), 
for each $\varepsilon \in (0,1]$ and $\lambda \in (0,1]$, 
there exist 
$\boldsymbol{v}_{\varepsilon ,\lambda } 
\in H^1(0,T;\boldsymbol{H}_0) \cap L^\infty (0,T;\mbox{\boldmath $ V$}_0)
\cap L^2(0,T;\boldsymbol{W})$ 
and $\mbox{\boldmath $ \mu $}_{\varepsilon ,\lambda } \in 
L^2(0,T;\mbox{\boldmath $ V$})$ 
such that
\begin{align} 
	\lambda \boldsymbol{v}'_{\varepsilon ,\lambda }(t) 
	& + \boldsymbol{F}^{-1} 
	\bigl( 
	\boldsymbol{v}'_{\varepsilon ,\lambda }(t)
	\bigr) 
	+ \varepsilon \partial \varphi \bigl( 
	\boldsymbol{v}_{\varepsilon ,\lambda } (t) \bigr) 
	\nonumber \\
	& 
	=\boldsymbol{P} \bigl( 
	-\boldsymbol{\beta }_\lambda \bigl( 
	\boldsymbol{u}_{\varepsilon ,\lambda }(t) \bigr)
	- \varepsilon \boldsymbol{\pi }
	\bigl( \boldsymbol{u}_{\varepsilon ,\lambda } (t) \bigr) +\boldsymbol{f}(t) \bigr)
	\quad {\rm in}~\boldsymbol{H}_0
	\label{ee}
\end{align} 
for a.a.\ $s\in (0,T)$
with $\boldsymbol{v}_{\varepsilon,\lambda } (0) 
=\boldsymbol{v}_{0\varepsilon }$ in 
$\boldsymbol{H}_0$, where 
$\boldsymbol{v}_{0\varepsilon }\in \boldsymbol{V}_0$ solves the auxiliary problem 
$\boldsymbol{v}_{0\varepsilon } + \varepsilon \partial 
\varphi (\boldsymbol{v}_{0\varepsilon })
=\boldsymbol{v}_0$ in $\boldsymbol{H}_0$ so that 
there exists a constant $C>0$ such that 
\begin{gather} 
	|\boldsymbol{v}_{0\varepsilon }|_{\boldsymbol{H}_0}^2 \le C,
	\quad 
	\varepsilon 
	|\boldsymbol{v}_{0\varepsilon }|_{\boldsymbol{V}_0}^2 \le C, 
	\label{constant}
	\\
	\int_{\Omega }^{} \widehat{\beta} (v_{0\varepsilon }+m_0) dx \le C,
	\quad 
	\int_{\Gamma }^{} \widehat{\beta} (v_{0\varepsilon }+m_0) d\Gamma  \le C.
	\nonumber 
\end{gather} 
Moreover, 
$\boldsymbol{u}_{\varepsilon, \lambda}
:=\boldsymbol{v}_{\varepsilon, \lambda} + m_0 \boldsymbol{1}$, 
$m_0:=m(\boldsymbol{u}_0)$ and 
$\boldsymbol{1}:=(1,1)$, 
and 
$\boldsymbol{\beta }_\lambda (\boldsymbol{z})
:= (\beta _\lambda (z),\beta _\lambda (z_\Gamma ))$ and 
$\boldsymbol{\pi }(\boldsymbol{z}):=(\pi (z),\pi (z_\Gamma ))$
for all $\boldsymbol{z} \in \mbox{\boldmath $ H $}$, where 
$\pi :D(\pi)=\mathbb{R} \to \mathbb{R}$ is 
a {L}ipschitz continuous function with a {L}ipschitz constant $L_\pi$ that breaks the monotonicity in  
$\beta + \varepsilon \pi $; 
$\boldsymbol{f} \in L^2(0,T;D(\partial \varphi ))$ is the 
solution of $\boldsymbol{g}(t)=\partial \varphi (\boldsymbol{f}(t))$ in 
$\boldsymbol{H}_0$ for a.a.\ $t \in (0,T)$. 
Namely, from \cite[Lemma~C]{CF15}, we can choose $\boldsymbol{f}(t):=(f(t),f_\Gamma (t))$ to satisfy
\begin{equation}
	\begin{cases}
	\displaystyle -\Delta f(t)=g(t) & {\rm a.e.\ in~} \Omega, \vspace{2mm}\\
	\displaystyle \partial _{\boldsymbol{\nu }}f(t)-\Delta _\Gamma f_\Gamma (t)=g_\Gamma (t) & {\rm a.e.\ on~} 
	\Gamma , 
	\end{cases} 
	\quad {\rm for~a.a.\ } t\in (0,T).
	\label{A3}
\end{equation}

\subsection{Uniform estimates for approximate solutions}

The key strategy in the proof is to obtain 
uniform estimates independent of $\varepsilon >0$ and $\lambda >0$, 
after which we consider the limiting procedures $\lambda \searrow 0$ and $\varepsilon \searrow 0$. 
Recall \eqref{ee} in the equivalent form
\begin{gather} 
	\boldsymbol{v}'_{\varepsilon ,\lambda }(s) 
	+ \boldsymbol{F} 
	\bigl( \boldsymbol{P} \boldsymbol{\mu}_{\varepsilon,\lambda  }(s)\bigr) 
	= \boldsymbol{0}
	\quad {\rm in}~ \boldsymbol{V}_0^*, 
	\label{ee2}\\
	\boldsymbol{\mu }_{\varepsilon ,\lambda }(s) 
	= \lambda \boldsymbol{v}_{\varepsilon, \lambda}'(s) 
	+ \varepsilon \partial  \varphi 
	\bigl( \boldsymbol{v}_{\varepsilon ,\lambda }(s) \bigr) 
	+
	\mbox{\boldmath $ \beta $}_\lambda \bigl( 
	\boldsymbol{u}_{\varepsilon ,\lambda }(s) 
	\bigr) 
	+ \varepsilon \boldsymbol{\pi }
	\bigl( \boldsymbol{u}_{\varepsilon ,\lambda } (s)
	\bigr)  -\boldsymbol{f}(s) 
	\quad {\rm in}~\boldsymbol{H}
	\label{ee3}
\end{gather} 
for a.a.\ $s \in (0,T)$.  
Moreover, if we put 
$\varepsilon _0:=\min \{1,1/(4L_\pi^2)\}$, then we have:

\paragraph{Lemma 3.1.}
{\it There exist 
positive constants $M_1$, $M_2$ 
independent of $\varepsilon \in (0,\varepsilon _0]$ and $\lambda \in (0,1]$ such that}
\begin{gather*}
	\lambda 
	\bigl| \boldsymbol{v}_{\varepsilon,\lambda}(t)
	\bigr|_{\boldsymbol{H}_0}^2 
	+ 
	\bigl| \boldsymbol{v}_{\varepsilon,\lambda }(t)
	\bigr|_{\boldsymbol{V}_0^*}^2
	\le M_1, \\
	\frac{\varepsilon}{2} \int_{0}^{t} 
	\bigl |\boldsymbol{v}_{\varepsilon,\lambda  }(s)
	\bigr |_{\boldsymbol{V}_0}^2 ds 
	+ 2 \int_{0}^{t} \bigl| \widehat{\beta}_\lambda 
	\bigl(u_{\varepsilon,\lambda }(s) \bigr) 
	\bigr|_{L^1(\Omega )} ds 
	+ 2 \int_{0}^{t} \bigl| 
	\widehat{\beta} _\lambda \bigl( u_{\Gamma, \varepsilon, \lambda}(s) \bigr) 
	\bigr|_{L^1(\Gamma)} ds
	\le M_2
\end{gather*}
{\it for all $t\in[0,T]$. }

\paragraph{Proof.} Multiplying \eqref{ee} by 
$\boldsymbol{v}_{\varepsilon,\lambda} (s) \in \boldsymbol{V}_0$, 
we have 
\begin{align*}
	& \lambda 
	\bigl( 
	\boldsymbol{v}_{\varepsilon ,\lambda }' (s),
	\boldsymbol{v}_{\varepsilon ,\lambda } (s) 
	\bigr)_{\boldsymbol{H}_0}  
	+  \bigl( 
	\boldsymbol{v}_{\varepsilon ,\lambda }' (s),
	\boldsymbol{v}_{\varepsilon ,\lambda } (s)
	\bigr )_{\boldsymbol{V}_0^*}
	+ \varepsilon \bigl( \partial \varphi \bigl ( \boldsymbol{v}_{\varepsilon,\lambda } (s) \bigr ), 
	\boldsymbol{v}_{\varepsilon ,\lambda }(s) 
	\bigr)_{\boldsymbol{H}_0} 
	\nonumber \\
	& \quad {} + \bigl( \boldsymbol{P} \boldsymbol{\beta }_\lambda
	\bigl ( \boldsymbol{v}_{\varepsilon,\lambda } (s)+m_0\boldsymbol{1} \bigr ), 
	\boldsymbol{v}_{\varepsilon,\lambda } (s) 
	\bigr)_{\boldsymbol{H}_0} 
	= 
	\bigl(\boldsymbol{f}(s)-\varepsilon \boldsymbol{P}
	\boldsymbol{\pi }
	\bigl( \boldsymbol{v}_{\varepsilon ,\lambda } (s)  +m_0 \boldsymbol{1}\bigr), 
	\boldsymbol{v}_{\varepsilon,\lambda } (s)
	\bigr)_{\boldsymbol{H}_0}
\end{align*} 
for a.a.\ $s \in (0,T)$. 
Using the definition of the subdifferential, 
we see that 
\begin{align*}
	& 
	\frac{\lambda }{2}\frac{d}{ds}
	\bigl | 
	\boldsymbol{v}_{\varepsilon, \lambda } (s )
	\bigr|_{\boldsymbol{H}_0}^2 
	+ 
	\frac{1}{2}
	\frac{d}{ds} 
	\bigl | 
	\boldsymbol{v}_{\varepsilon, \lambda }(s )
	\bigr |_{\boldsymbol{V}_0^*}^2
	+ \frac{\varepsilon}{2} 
	\bigl | \boldsymbol{v}_{\varepsilon,\lambda }(s )
	\bigr |_{\boldsymbol{V}_0}^2
	\nonumber \\
	&	
	\quad {} +
	\bigl| \widehat{\beta}_\lambda 
	\bigl ( u_{\varepsilon,\lambda } (s) 
	\bigr) \bigr |_{L^1(\Omega )}
	+ \bigl| \widehat{\beta}_\lambda 
	\bigl ( u_{\Gamma, \varepsilon,\lambda } (s) 
	\bigr) \bigr |_{L^1(\Gamma )}
	\nonumber 
	\\
	& \le 
	\bigl( |\Omega |+|\Gamma | \bigr) \widehat{\beta }(m_0)
	+
	\frac{1}{2} \bigl| \boldsymbol{v}_{\varepsilon ,\lambda } (s) 
	\bigr|_{\boldsymbol{V}_0^*}^2  
	+ 
	\bigl| \boldsymbol{f}(s) \bigr|_{\boldsymbol{V}_0}^2
	+ L_\pi ^2 \varepsilon ^2 \bigl| 
	\boldsymbol{v}_{\varepsilon,\lambda }(s)
	\bigr|_{\boldsymbol{V}_0}^2
	\quad {\rm for~a.a.\ } s \in (0,T).
\end{align*}
Taking $\varepsilon \in (0,\varepsilon _0]$ and 
using the Gronwall inequality, we obtain 
the existence of $M_1$ and $M_2$ 
independent of $\varepsilon \in (0,\varepsilon _0]$ and $\lambda \in (0,1]$
satisfying the conclusion. \hfill $\Box$

\paragraph{Lemma 3.2.}
{\it There exists a positive constant $M_3$, 
independent of $\varepsilon \in (0,\varepsilon _0]$ and $\lambda \in (0,1]$, such that}
\begin{align*}
	& 
	2 \lambda 
	\int_{0}^{t} 
	\bigl| 
	\boldsymbol{v}_{\varepsilon,\lambda } '(s) 
	\bigr|_{\boldsymbol{H}_0}^2 
	ds
	+ 
	\int_{0}^{t} 
	\bigl| 
	\boldsymbol{v}_{\varepsilon,\lambda } '(s)
	\bigr| _{\boldsymbol{V}_0^*}^2 ds 
	+ 
	\varepsilon 
	\bigl|\boldsymbol{v}_{\varepsilon,\lambda }(t) 
	\bigr|_{\boldsymbol{V}_0}^2  
	\nonumber
	\\ 
	&{}+ 2 \bigl| \widehat{\beta}_\lambda  \bigl (u_{\varepsilon,\lambda}(t) \bigr ) 
	\bigr|_{L^1(\Omega )} 
	+ 2 \bigl| \widehat{\beta}_\lambda  \bigl ( 
	u_{\Gamma, \varepsilon,\lambda }(t) \bigr) 
	\bigr|_{L^1(\Gamma )} \le M_3, 
\end{align*}
\begin{equation*} 
	 \int_{0}^{t} \bigl| \boldsymbol{P}
	 \boldsymbol{\mu }_{\varepsilon, \lambda }(s) 
	 \bigr|_{\boldsymbol{V}_0}^2 ds 
	 \le M_3 \quad {\it for~all~} t\in [0,T].
\end{equation*}

\paragraph{Proof.} 
Multiplying \eqref{ee} by 
$\boldsymbol{v}_{\varepsilon,\lambda }' (s) 
\in \boldsymbol{H}_0$, we have 
\begin{align*}
	& \lambda \bigl 
	|\boldsymbol{v}_{\varepsilon,\lambda } '(s ) 
	\bigr |_{\boldsymbol{H}_0}^2  
	+  
	\frac{1}{2} \bigl |\boldsymbol{v}_{\varepsilon,\lambda } '(s) 
	\bigr |_{\boldsymbol{V}_0^*}^2 
	+ \varepsilon \frac{d}{ds} 
	\varphi \bigl ( \boldsymbol{v}_{\varepsilon,\lambda } (s) \bigr )
	+\frac{d}{ds} \int_{\Omega }^{} 
	\widehat{\beta }_\lambda  \bigl (
	u_{\varepsilon,\lambda } (s) \bigr)dx
	\nonumber \\
	&\quad {}
	+
	\frac{d}{ds} \int_{\Gamma }^{} 
	\widehat{\beta }_\lambda  
	\bigl (u_{\Gamma, \varepsilon,\lambda } (s)\bigr)d\Gamma 
	\le 
	L_\pi ^2 \varepsilon ^2 \bigl| \boldsymbol{v}_{\varepsilon ,\lambda }(s) 
	\bigr|_{\boldsymbol{V}_0}^2
	+ \bigl| \boldsymbol{f}(s) \bigr|_{\boldsymbol{V}_0}^2  
	\quad {\rm for~a.a.\ } s\in (0,T).
\end{align*}
Integrating this over $(0,t)$ with respect to $s$, 
we see that there exists a positive constant $M_3$, independent 
of $\varepsilon \in (0,\varepsilon _0]$ and $\lambda \in (0,1]$, such that 
the first estimate holds.  
Next, multiplying \eqref{ee2} by 
$\boldsymbol{P} \boldsymbol{\mu }_{\varepsilon,\lambda }(s) 
\in \boldsymbol{V}_0$ and integrating the 
resultant over $(0,t)$ with respect to $s$, 
we obtain the second estimate.
\hfill $\Box$ \\

The previous two lemmas are essentially the same as \cite[Lemmas~3.1 and 3.2]{Fuk15}. 
The next uniform estimate is the point of emphasis in this paper.

\paragraph{Lemma 3.3.}
{\it There exists positive constant $M_4$,
independent of $\varepsilon \in (0,1]$ and $\lambda \in (0,1]$, such that}
\begin{gather*}
	 \bigl| 
	 \boldsymbol{u}_{\varepsilon, \lambda }(t) 
	 \bigr|_{\boldsymbol{H}}^2
	 \le M_4 \left( 1+ \frac{\lambda }{\varepsilon }\right),
	 \quad \bigl| 
	 \boldsymbol{v}_{\varepsilon, \lambda }(t) 
	 \bigr|_{\boldsymbol{H}_0}^2
	 \le M_4\left( 1+ \frac{\lambda }{\varepsilon }\right),
	 \\
	 \lambda 
	 \bigl| \boldsymbol{v}_{\varepsilon ,\lambda }(t) 
	 \bigr|_{\boldsymbol{V}_0}^2 
	 + \varepsilon \int_{0}^{t} 
	 \bigl|
	 \partial \varphi 
	 \bigl( 
	 \boldsymbol{v}_{\varepsilon ,\lambda }(s)
	 \bigr) 
	 \bigr|^2_{\boldsymbol{H}_0} ds 
	 \le M_4\left( 1+ \frac{\lambda }{\varepsilon }\right)
	 \quad {\it for~all~} t\in [0,T].
\end{gather*}

\paragraph{Proof.} Multiplying \eqref{ee2} by 
$\boldsymbol{v}_{\varepsilon,\lambda}(s) \in \boldsymbol{V}_0$ and 
using the fact $(d/ds)m(\boldsymbol{u}_{\varepsilon,\lambda }(s))=0$, 
 we have 
\begin{equation*} 
	\bigl( \boldsymbol{u}'_{\varepsilon,\lambda }(s), 
	\boldsymbol{u}_{\varepsilon,\lambda } (s)
	\bigr)_{\boldsymbol{H}}
	+ a \bigl(\boldsymbol{\mu}_{\varepsilon,\lambda  }(s), 
	\boldsymbol{u}_{\varepsilon ,\lambda }(s) \bigr)
	=0
\end{equation*} 
for a.a.\ $s \in (0,T)$ 
(see \cite[Remark~3]{Fuk15}). 
On the other hand, multiplying \eqref{ee3} by 
$\partial \varphi (\boldsymbol{v}_{\varepsilon,\lambda}(s)) \in \boldsymbol{H}_0$ 
and integrating by parts, we have
\begin{align*}
	a \bigl(\boldsymbol{\mu}_{\varepsilon,\lambda  }(s), 
	\boldsymbol{u}_{\varepsilon ,\lambda }(s) \bigr) 
	& = 
	\frac{\lambda }{2} \frac{d}{ds} a 
	\bigl(\boldsymbol{u}_{\varepsilon,\lambda  }(s), 
	\boldsymbol{u}_{\varepsilon ,\lambda }(s) \bigr) 
	+  \varepsilon 
	\bigl| 
	\partial \varphi 
	\bigl(
	\boldsymbol{v}_{\varepsilon,\lambda  }(s)
	\bigr)
	 \bigr| _{\boldsymbol{H}_0}^2
	\nonumber \\
	& \quad {} + \int_{\Omega }^{} \beta '_\lambda 
	 \bigl( 
	 u_{\varepsilon ,\lambda }(s) 
	 \bigr) 
	 \bigl| 
	 \nabla u_{\varepsilon ,\lambda }(s)
	 \bigr| ^2 dx 
	 +
	\int_{\Gamma }^{} \beta '_\lambda 
	 \bigl( 
	 u_{\Gamma, \varepsilon ,\lambda }(s) 
	 \bigr) 
	 \bigl| 
	 \nabla_\Gamma  u_{\Gamma, \varepsilon ,\lambda }(s)
	 \bigr| ^2 d\Gamma
	 \nonumber \\
	 & \quad {}
	 + \varepsilon 
	 \bigl( \boldsymbol{\pi } 
	 \bigl( \boldsymbol{u}_{\varepsilon,\lambda  }(s)
	 \bigr), 
	 \partial \varphi 
	 \bigl( \boldsymbol{u}_{\varepsilon ,\lambda }(s) \bigr) 
	 \bigr)_{\boldsymbol{H}}
	 +
	\bigl(\partial \varphi 
	\bigl( 
	\boldsymbol{f}(s)
	\bigr),
	\boldsymbol{u}_{\varepsilon ,\lambda }(s)
	\bigr)_{\mbox{\scriptsize \boldmath $ H$}}
\end{align*} 
for a.a.\ $s \in (0,T)$. 
Using the {L}ipschitz continuity of $\pi $ and \eqref{A3}, we see that there exists a positive 
constant $C_\pi $ such that
\begin{equation*}
	\frac{d}{ds} 
	\bigl| 
	\boldsymbol{u}_{\varepsilon,\lambda }(s)
	\bigr|_{\boldsymbol{H}}^2 
	+ 
	\lambda 
	\frac{d}{ds} 
	\bigl| 
	\boldsymbol{v}_{\varepsilon,\lambda }(s)
	\bigr|_{\boldsymbol{V}_0}^2 
	+ \varepsilon 
	\bigl| 
	\partial \varphi 
	\bigl(
	\boldsymbol{v}_{\varepsilon,\lambda  }(s)
	\bigr)
	 \bigr| _{\boldsymbol{H}_0}^2
	 \le C_\pi \bigl( \bigl| 
	 \boldsymbol{u}_{\varepsilon ,\lambda }(s) \bigr|_{\boldsymbol{H}}^2 
	 + 1 \bigr) 
	 + \bigl| 
	 \boldsymbol{g}(s)
	 \bigr|_{\boldsymbol{H}_0}^2
\end{equation*} 
for a.a.\ $s \in (0,T)$. 
Then, using \eqref{constant} and the {G}ronwall inequality, we deduce that
\begin{align*} 
	\bigl| 
	\boldsymbol{u}_{\varepsilon,\lambda }(t)
	\bigr|_{\boldsymbol{H}}^2 
	+ 
	\lambda 
	\bigl| 
	\boldsymbol{v}_{\varepsilon,\lambda }(t)
	\bigr|_{\boldsymbol{V}_0}^2 
	& \le 
	\left\{ 
	| 
	\boldsymbol{v}_{0\varepsilon}+m_0 \boldsymbol{1}
	|_{\boldsymbol{H}}^2 
	+ 
	\lambda 
	| 
	\boldsymbol{v}_{0\varepsilon }
	|_{\boldsymbol{V}_0}^2 
	 + C_\pi  T
	 + | 
	 \boldsymbol{g}
	 |_{L^2(0,T;\boldsymbol{H}_0)}^2
	\right\}
	e^{C_\pi T} 
	\\
	& \le \left\{ 
	2C+2|m_0|^2 \bigl( |\Omega |+|\Gamma | \bigr) 
	+ 
	\frac{\lambda }{\varepsilon } C
	 + C_\pi  T
	 + | 
	 \boldsymbol{g}
	 |_{L^2(0,T;\boldsymbol{H}_0)}^2
	\right\}
	e^{C_\pi T}
\end{align*} 
for all $t\in [0,T]$. That is, there exists a positive constant $M_4$ independent of 
$\varepsilon \in (0,1]$ and $\lambda \in (0,1]$ such that
the uniform estimates hold. 
\hfill $\Box$ 

\paragraph{Lemma 3.4.}
{\it There exists positive constant $M_5$,
independent of $\varepsilon \in (0,1]$ and $\lambda \in (0,1]$, such that}
\begin{gather*}
	\int_{0}^{t}
	\bigl| \boldsymbol{\mu }_{\varepsilon,\lambda }(s) \bigr|_{\boldsymbol{V}_0}^2 ds 
	\le M_5 \left( 1+ \frac{\lambda }{\varepsilon }\right), 
	\\
	\int_{0}^{t} 
	\bigl| \boldsymbol{\beta }_\lambda \bigl 
	(\boldsymbol{u}_{\varepsilon,\lambda }(s) \bigr) 
	\bigr|_{\boldsymbol{H}}^2 ds 
	\le M_5 \left( 1+ \frac{\lambda }{\varepsilon }\right)
	\quad {\it for~all~} t\in [0,T].
\end{gather*}

Using Lemmas 3.1 to 3.3,  
the proofs of these uniform estimates are completely the same as those for 
\cite[Lemmas~4.3 and 4.4]{CF15}. 
Therefore, we omit the proof.

\subsection{Limiting procedure}

From the previous uniform estimates, we can 
consider the limit as $\lambda \searrow 0$. More precisely, 
for each $\varepsilon \in (0,\varepsilon _0]$,
there exists a subsequence $\{ \lambda_k \}_{k \in \mathbb{N}}$ with 
$\lambda _k \searrow 0$ as $k \to +\infty $ and 
a quadruplet $(\boldsymbol{v}_\varepsilon, 
\boldsymbol{v}_\varepsilon ^*,
\boldsymbol{\mu }_\varepsilon, \boldsymbol{\xi }_\varepsilon )$ 
of $\boldsymbol{v}_\varepsilon 
\in H^1(0,T;\boldsymbol{V}_0^*) 
\cap L^\infty (0,T;\boldsymbol{V}_0) 
\cap L^2(0,T;\boldsymbol{W})$, 
$\boldsymbol{v}_\varepsilon ^*
\in L^2(0,T;\boldsymbol{H}_0)$, 
$\boldsymbol{\mu }_\varepsilon 
\in L^2(0,T;\boldsymbol{V})$, 
$\boldsymbol{\xi }_\varepsilon 
\in L^2(0,T;\boldsymbol{H})$, 
such that 
\begin{gather*} 
	\boldsymbol{v}_{\varepsilon, \lambda_k}
	\to \boldsymbol{v}_\varepsilon 
	\quad \mbox{weakly star in } 
	L^\infty (0,T;\boldsymbol{H}_0), 
	\quad 
	\boldsymbol{v}_{\varepsilon, \lambda _k} 
	\to \boldsymbol{v} _\varepsilon 
	\quad \mbox{weakly in } 
	L^2 (0,T;\boldsymbol{V}_0), 
	\\
	\lambda _k \boldsymbol{v}_{\varepsilon, \lambda _k}' 
	\to \mbox{\boldmath $ 0$}
	\quad \mbox{in } 
	L^2 (0,T;\boldsymbol{H}_0), 
	\quad 
	\boldsymbol{v}_{\varepsilon, \lambda _k}' \to \boldsymbol{v}_\varepsilon ' 
	\quad \mbox{weakly in } 
	L^2(0,T;\boldsymbol{V}_0^*),
	\\
	\boldsymbol{u}_{\varepsilon, \lambda _k} \to \boldsymbol{u}_\varepsilon :=
	\boldsymbol{v}_\varepsilon +m_0\boldsymbol{1}
	\quad \mbox{weakly star in } 
	L^\infty (0,T;\boldsymbol{V}),
	\\
	\partial \varphi  (\boldsymbol{v}_{\varepsilon ,\lambda _k}) 
	\to \boldsymbol{v}_\varepsilon ^* 
	\quad \mbox{weakly in } 
	L^2 (0,T;\boldsymbol{H}_0),
	\quad 
	\boldsymbol{\mu }_{\varepsilon, \lambda _k} \to 
	\boldsymbol{\mu } _\varepsilon 
	\quad \mbox{weakly in } 
	L^2(0,T;\boldsymbol{V}), 
	\\
	\boldsymbol{\beta }_{\lambda _k} 
	( 
	\boldsymbol{u}_{\varepsilon, \lambda _k} 
	) 
	\to \boldsymbol{\xi }_\varepsilon  \quad \mbox{weakly in } 
	L^2(0,T;\boldsymbol{H})
	\quad {\rm as~} k \to +\infty.
\end{gather*} 
From the compactness theorem
(see, e.g., \cite[Section~8, Corollary~4]{Sim87}), this 
gives
\begin{gather*} 
	\boldsymbol{v}_{\varepsilon, \lambda _k} 
	\to \boldsymbol{v}_\varepsilon  \quad \mbox{in } 
	C\bigl( [0,T];\boldsymbol{H}_0 \bigr),
	\quad 
	\boldsymbol{u}_{\varepsilon, \lambda _k} 
	\to \boldsymbol{u}_\varepsilon  \quad \mbox{in } 
	C\bigl( [0,T];\boldsymbol{H} \bigr),
	\\
	\mbox{\boldmath $\pi $} 
	(\boldsymbol{u}_{\varepsilon, \lambda _k}) 
	\to \mbox{\boldmath $\pi $} 
	(\boldsymbol{u}_{\varepsilon}) \quad \mbox{in } 
	C\bigl( [0,T];\boldsymbol{H} \bigr)
	\quad {\rm as~} k \to +\infty.
\end{gather*}
Moreover, 
from the demi-closedness of $\partial \varphi $ and 
\cite[Proposition~2.2]{Bar10}, 
we see that 
$\boldsymbol{v}_\varepsilon^* =\partial \varphi (\boldsymbol{v}_\varepsilon )$ 
in $L^2(0,T;\boldsymbol{H}_0)$
and 
$\boldsymbol{\xi }_\varepsilon \in \boldsymbol{\beta }(\boldsymbol{u}_\varepsilon )$ 
in $L^2(0,T;\boldsymbol{H})$. From these facts, we deduce from \eqref{ee2} and \eqref{ee3} that
\begin{gather} 
	\boldsymbol{v}'_{\varepsilon}(t) 
	+ \boldsymbol{F} 
	\bigl( \boldsymbol{P} \boldsymbol{\mu}_{\varepsilon}(t)\bigr) 
	= \boldsymbol{0}
	\quad {\rm in}~ \boldsymbol{V}_0^*, 
	\label{ee4}\\
	\boldsymbol{\xi }_\varepsilon (t) \in 
	\boldsymbol{\beta }\bigl(\boldsymbol{u}_\varepsilon (t)\bigr), 
	\quad 
	\boldsymbol{\mu }_{\varepsilon}(t) 
	= 
	\varepsilon \partial  \varphi 
	\bigl( \boldsymbol{v}_{\varepsilon}(t) \bigr) 
	+
	\boldsymbol{\xi }_\varepsilon (t) 
	+ \varepsilon \boldsymbol{\pi }
	\bigl( \boldsymbol{u}_{\varepsilon} (t)
	\bigr)  -\boldsymbol{f}(t) 
	\quad {\rm in}~\boldsymbol{H}
	\label{ee5}
\end{gather} 
for a.a.\ $t \in (0,T)$, with 
$\boldsymbol{v}_\varepsilon (0) = \boldsymbol{v}_{0\varepsilon }$ 
in $\boldsymbol{H}$. 
We also have the regularity 
$\boldsymbol{u}_\varepsilon 
\in H^1(0,T;\boldsymbol{V}^*) \cap L^\infty (0,T;\boldsymbol{V}) 
\cap L^2(0,T;\boldsymbol{W})$. 
Now, taking 
the limit inferior as $\lambda \searrow 0$ on the uniform estimates, $\lambda /\varepsilon \searrow 0$ for 
all $\varepsilon \in (0,\varepsilon _0]$, and 
we therefore obtain the same kind of uniform estimates 
as in the previous lemmas independent of $\varepsilon \in (0,\varepsilon _0]$.

\paragraph{Proof of Theorem 2.1.}
By using the estimates for 
$\boldsymbol{v}_\varepsilon, 
\boldsymbol{u}_\varepsilon, 
\boldsymbol{\mu }_\varepsilon$ 
and 
$\boldsymbol{\xi}_\varepsilon$, 
there 
exist a subsequence $\{ \varepsilon _k \}_{k \in \mathbb{N}}$ with 
$\varepsilon _k \searrow 0$ as $k \to +\infty $
and functions 
$\boldsymbol{v} \in H^1(0,T;\boldsymbol{V}_0^*) \cap L^\infty (0,T;\boldsymbol{H}_0)$, 
$\boldsymbol{u} \in H^1(0,T;\boldsymbol{V}^*) 
\cap L^\infty (0,T;\boldsymbol{H})$, 
$\mbox{\boldmath $ \mu $}\in L^2(0,T;\boldsymbol{V})$ and 
$\mbox{\boldmath $ \xi $} \in L^2(0,T;\boldsymbol{H})$ such that 
\begin{gather*} 
	\boldsymbol{v}_{\varepsilon_k} 
	\to \boldsymbol{v} \quad \mbox{weakly star in } 
	L^\infty (0,T;\boldsymbol{H}_0), 
	\\
	\boldsymbol{u}_{\varepsilon_k} 
	\to \boldsymbol{u}=\boldsymbol{v}+m_0 \mbox{\boldmath $ 1$} \quad \mbox{weakly star in } 
	L^\infty (0,T;\boldsymbol{H}), 
	\quad 
	\varepsilon _k \boldsymbol{v}_{\varepsilon_k} \to \mbox{\boldmath $ 0$}
	\quad \mbox{in } 
	L^\infty (0,T;\boldsymbol{V}_0), 
	\\
	\boldsymbol{v}_{\varepsilon_k}' \to \boldsymbol{v}' 
	\quad \mbox{weakly in } 
	L^2(0,T;\boldsymbol{V}_0^*),
	\quad 
	\boldsymbol{u}_{\varepsilon_k}' \to \boldsymbol{u}' 
	\quad \mbox{weakly in } 
	L^2(0,T;\boldsymbol{V}^*),
	\\
	\boldsymbol{\mu }_{\varepsilon_k} \to \boldsymbol{\mu }
	\quad \mbox{weakly in } 
	L^2(0,T;\boldsymbol{V}), 
	\quad 
	\boldsymbol{\xi }_{\varepsilon_k} 
	\to \boldsymbol{\xi } \quad \mbox{weakly in } 
	L^2(0,T;\boldsymbol{H}),
	\\
	\varepsilon_k \boldsymbol{\pi } 
	(\boldsymbol{u}_{\varepsilon_k}) 
	\to \boldsymbol{0} \quad \mbox{in } 
	L^\infty (0,T;\boldsymbol{H}) 
	\quad {\rm as~} k \to +\infty.
\end{gather*} 
From 
the {A}scoli--{A}rzel\`{a} theorem, we also have
\begin{gather*} 
	\boldsymbol{v}_{\varepsilon_k} \to \boldsymbol{v} \quad \mbox{in } 
	C\bigl( [0,T];\boldsymbol{V}_0^* \bigr),
	\quad 
	\boldsymbol{u}_{\varepsilon_k} \to \boldsymbol{u} \quad \mbox{in } 
	C\bigl( [0,T];\boldsymbol{V}^* \bigr)
	\quad {\rm as~} k \to +\infty.
\end{gather*}
Now, multiplying \eqref{ee5} by 
$\boldsymbol{\eta } \in L^2(0,T;\boldsymbol{V})$ and 
integrating over $(0,T)$, we obtain
\begin{align} 
	\int_{0}^{T} \bigl( 
	\boldsymbol{\mu }_{\varepsilon_k} (t),\boldsymbol{\eta}(t)
	\bigr)_{\boldsymbol{H}} dt
	& = \varepsilon_k \int_{0}^{T} 
	a \bigl(  \boldsymbol{v}_{\varepsilon_k} (t), 
	\boldsymbol{\eta }(t) \bigr) dt
	+ \int_{0}^{T} 
	\bigl( 
	 \boldsymbol{\xi }_{\varepsilon_k} (t), 
	 \boldsymbol{\eta }(t) 
	 \bigr) _{\boldsymbol{H}}
	 dt 
	 \nonumber \\
	 & \quad {}
	+ \varepsilon_k \int_{0}^{T}
	\bigl( 
	\boldsymbol{\pi }\bigl( 
	\boldsymbol{u}_{\varepsilon_k} (t)\bigr), \boldsymbol{\eta }(t)
	\bigr)_{\boldsymbol{H}} 
	dt
	-
	\int_{0}^{T}
	\bigl( \boldsymbol{f}(t),\boldsymbol{\eta }(t)
	\bigr)_{\boldsymbol{H}} dt.
	\label{weakvari}
\end{align}
Letting $k \to \infty $, we obtain 
\begin{equation*} 
	\int_{0}^{T} \bigl( 
	\boldsymbol{\mu } (t),\boldsymbol{\eta }(t)
	\bigr)_{\boldsymbol{H}} dt
	= 
	\int_{0}^{T} 
	\bigl( 
	 \boldsymbol{\xi }(t)-\boldsymbol{f}(t),\boldsymbol{\eta }(t)
	\bigr)_{\boldsymbol{H}} dt
	\quad \mbox{\rm for all }\boldsymbol{\eta } \in L^2(0,T;\boldsymbol{V}),
\end{equation*}
namely, $\boldsymbol{\mu }=\boldsymbol{\xi }-\boldsymbol{f}$ in 
$L^2(0,T;\boldsymbol{H})$. 
This implies 
the regularity of $\boldsymbol{\xi } \in L^2(0,T;\boldsymbol{V})$,
that is, $\xi _\Gamma =\xi _{|_\Gamma }$ a.e.\ on $\Sigma $. 
Next, we take $\boldsymbol{\eta}:=\boldsymbol{u}_{\varepsilon _k} 
\in L^2(0,T;\boldsymbol{V})$ in 
\eqref{weakvari}, so that 
\begin{align*} 
	\limsup _{k \to +\infty }
	\int_{0}^{T} 
	\bigl( 
	 \boldsymbol{\xi }_{\varepsilon_k} (t), 
	 \boldsymbol{u}_{\varepsilon _k}(t) 
	 \bigr) _{\boldsymbol{H}}
	 dt 
	 & \le 
	\int_{0}^{T} \bigl \langle \boldsymbol{u}(t),
	\boldsymbol{\mu }(t)
	\bigr \rangle_{\boldsymbol{V}^*, 
	\boldsymbol{V}} dt
	+ 
	\int_{0}^{T}
	\bigl( \boldsymbol{f}(t),\boldsymbol{u}(t)
	\bigr)_{\boldsymbol{H}} dt \\
	& = \int_{0}^{T} \bigl( 
	 \boldsymbol{\xi }(t), 
	 \boldsymbol{u}(t) 
	 \bigr) _{\boldsymbol{H}}
	 dt. 
\end{align*}
Thus, applying \cite[Proposition~2.2]{Bar10} we have 
$\boldsymbol{\xi } \in \boldsymbol{\beta }(\boldsymbol{u})$ 
in $L^2(0,T;\boldsymbol{H})$, and so 
we obtain 
$\xi \in \beta (u)$ a.e.\ in $Q$.
$\xi _\Gamma \in \beta (u_\Gamma )$ a.e.\ on $\Sigma$.
Finally, letting $k \to +\infty $ and applying 
{H}ahn--{B}anach extension theorem of bounded linear functional on 
$\boldsymbol{V}$ to $V \times V_\Gamma$,
then we see that \eqref{ee4} gives \eqref{def}
for a.a.\ $t \in (0,T)$, with $u(0)=u_0$ 
a.e.\ in $\Omega $ and $u_\Gamma (0)=u_{0\Gamma }$ a.e.\ on $\Gamma $. 
\hfill $\Box$


\section*{Acknowledgments}

The author is indebted to professor 
Pierluigi Colli, who kindly gave him the opportunity for fruitful discussions. 
The author is supported by JSPS KAKENHI 
Grant-in-Aid for Scientific Research(C), Grant Number 26400164.  
Last but not least, the author is also grateful to the referee for the careful reading of the
manuscript.


\end{document}